\documentclass[12pt]{amsart}
\usepackage{amscd,amsthm,amssymb,amsfonts,amsmath}
\usepackage{bbm}
\usepackage{url}

\theoremstyle{plain}
\newtheorem{thm}{Theorem}[section]
\newtheorem{lemma}[thm]{Lemma}
\newtheorem{prop}[thm]{Proposition}
\newtheorem{cor}[thm]{Corollary}
\newtheorem{conjecture}[thm]{Conjecture}
\theoremstyle{definition}
\newtheorem*{defn}{Definition}

\theoremstyle{remark}
\newtheorem*{remark}{Remark}
\newtheorem*{ack}{Acknowledgments}
\newtheorem*{notations}{Notations}

\newcommand{\Z}{\mathbb Z}    
\newcommand{\R}{\mathbb R}    
\newcommand{\C}{\mathbb C}    
\newcommand{\Hess}{\operatorname{Hess}} 

\newcommand{\<}{\langle}   
\renewcommand{\>}{\rangle} 
\newcommand{\suchthat}{\ : \ }

\newcommand{\nc}{\Pi}
\newcommand{\cone}{\operatorname{cone}}

\newcommand{\dU}{\partial\mathcal U}

\begin{document} 

\begin{abstract}
We use a generalization of the Gibbons-Hawking ansatz to study the behavior of certain non-compact Calabi-Yau manifolds in the large complex structure limit. This analysis provides an intermediate step toward proving the metric collapse conjecture for toric hypersurfaces and complete intersections.  
\end{abstract}
 
\title[Limiting behavior of local Calabi-Yau metrics]{Limiting behavior of local Calabi-Yau  metrics
}
\author{Ilia Zharkov}
\address{Mathematics Department \\ Harvard University \\ Cambridge, MA 02138
  \\ USA}
\email{zharkov@math.harvard.edu}
\maketitle

\section{Introduction} 
Since the Strominger-Yau-Zaslow conjecture \cite{SYZ} was made there has been a considerable interest in geometry of K\"ahler Ricci-flat half-dimensional torus fibrations. For a nonsingular fibration the metrics which are flat along the fibers were used by Hitchin \cite{Hi} in relation with mirror symmetry and SYZ conjecture. This, so called, semi-flat case was further studied in \cite{LYZ} and \cite{Leung}. The Ricci-flatness condition becomes equivalent to the real Monge-Amp\`ere equation on the base. The mirror symmetry is then provided by the Legendre transform.

This simple case has inspired another, more recent, conjecture of Gross and Wilson \cite{GW} and Kontsevich and Soibelman \cite{KS} about existence of an integral K\"ahler affine structure on the limiting space of the metric collapse. But in order to understand this conjecture for general compact Calabi-Yau manifolds (not just tori) one needs a local description of the metric behavior near the singular fibers. 

The first example of such a metric was constructed by Ooguri and Vafa \cite{OV} in two dimensions via periodic Gibbons-Hawking ansatz and was used later by Gross and Wilson \cite{GW} to justify the collapsing of K3 to $S^2$ with explicit metric (singular at 24 points). The limiting metric on the $S^2$ was, in fact, written down long before the SYZ story by Greene et al. \cite{stringy}. In higher dimensions Pedersen and Poon \cite{PP} derived the (non-linear) differential equations for the GH ansatz but the solutions they found were not periodic in the remaining number of torus variables, and hence cannot apply to our case of interest. An attempt to apply the generalization of Gibbons-Hawking in the periodic situation in dimension 3 was made by Matessi \cite{Matessi}, though no explicit solutions were found.

In this paper we do not try to solve Gibbons-Hawking differential equation. Rather the goal is to set up the geometric framework  for investigation of limiting behavior of such metrics as the size of tori goes to zero. Unfortunately, the key exponential decay lemma is left unproven. A substantial amount of hard analysis of non-linear elliptic PDE with singularities is required for the proof, and we plan to do it elsewhere.

The conjectural description of the limiting metric space is in agreement with the metric collapse picture of \cite{KS} and \cite{GW}. But we suggest a more explicit condition for asymptotics at the singular locus. The coordinates of the Gibbons-Hawking ansatz are related to the affine ones via the partial Legendre transform introduced in the last section of the paper.

\begin{notations}
We use the usual convention to sum over repeated indices. Also we deal with orbifolds on the same footing as with regular complex manifolds. That is, when we say a form or a map is holomorphic, it is meant in this orbifold sense. 
\end{notations}

\begin{ack}
I am very grateful to M. Gross for explaining to me the right normalization of the holomorphic volume form and several other key issues. I have also considerably benefited from conversations with R. Bryant, C. Haase, M. Kontsevich, D. Morrison, M. Stern and S.-T. Yau. I am indebted as well to the referee for a whole series of important comments and corrections to the original text. 
Finally, I would like to thank IHES for its hospitality and financial support during  the summers of 2001 and 2002 where a significant part of the work has been done.
\end{ack}

\section{Ricci-flat metrics} 

\subsection{Generalized Gibbons-Hawking ansatz}
Suppose $T^n$, a $n$-dimensional (real) torus, acts freely on an $N$-(complex) dimensional K\"ahler manifold $M$ by Hamiltonian holomorphic isometries. Then $M$ can be considered as a principal $T^n$-bundle over a real manifold of dimension $2N-n$. The generalized Gibbons-Hawking ansatz expresses the K\"ahler and Ricci-flat conditions as differential equations in the $n$ moment map coordinates and $N-n$ holomorphic coordinates on the K\"ahler quotient. 

Let $\mathfrak t$ denote the Lie algebra of the torus Lie group $T^n$, and let $\mathfrak t_\Z$ be the natural integral lattice in $\mathfrak t$. We will fix a basis in $\mathfrak t_\Z$. This defines affine coordinates $u_i$ on the dual space $\mathfrak t^*\cong \R^n$.
Let $Y$ be either $\C^{N-n}$ or $(\C^*)^{N-n}\cong\R^{N-n}\times(S^1)^{N-n}$, with the affine complex coordinates $\eta_p=x_p+iy_p$, where $y_p$ are the phase coordinates on the torus $(S^1)^{N-n}$ in the latter case.

Consider a principal $T^n$-bundle  $\pi:M\to B^\circ$ over an open set $B^\circ$ in $\mathfrak t^* \times Y$, with coordinates $(u_i,\eta_p,\bar{\eta}_q)$. Denote by $[\nu]$ its integral Chern class as an element in $H^2 (B^\circ,\mathfrak t_\Z)$.

\begin{thm}[cf. \cite{PP}]\label{thm:GH}
Let $V^{ij}$, respectively  $W^{pq}$, be real symmetric, respectively hermitian,  positive definite matrices of smooth functions on $B^\circ$, locally given by some potential function $\Phi$:
 \begin{equation} \label{eq:potential}
   V^{ij}=\frac{\partial^2\Phi}{\partial u_i\partial u_j}, \quad        
   W^{pq}=-4\frac{\partial^2\Phi}{\partial \eta_p \partial \bar{\eta}_q}, \qquad 1\le i,j \le n \quad n+1\le p,q \le N.
 \end{equation}
Then the following $\mathfrak t$-valued 2-form is closed:
 \begin{equation}
   F_j=\sqrt{-1}\left(\frac12 \frac{\partial W^{pq}}{\partial u_j}d\eta_p\wedge
   d\bar{\eta}_q +\frac{\partial V^{ij}}{\partial \eta_p}du_i\wedge d\eta_p 
   -\frac{\partial V^{ij}}{\partial \bar{\eta}_q}du_i\wedge d\bar{\eta}_q 
   \right).
\end{equation}
Suppose, in addition, that $\det V^{ij}=\det W^{pq}$ and $(F_1,\dots,F_n)$ is in the cohomology class $2\pi[\nu]$. Then there exist a connection on the bundle $M\to B^\circ$ with associated 1-forms $A_i$ and the curvature $(F_1,\dots,F_n)$ such that $M$ is a K\"ahler manifold with Ricci-flat metric given by
 \begin{equation}
   h=(V^{-1})^{ij}dz_i\otimes d\bar{z}_j+W^{pq} d\eta_p\otimes d\bar{\eta}_q,
 \end{equation}    
where $dz_j=V^{ij}du_i+\sqrt{-1}\cdot A_j$ and $d\eta_p$ form a basis of holomorphic 1-forms. The holomorphic $N$-form and the K\"ahler form:
\begin{equation}
  \Omega=\wedge_{j=1}^k dz_j\bigwedge \wedge_{p=1}^l d\eta_p , \qquad 
  \omega=du_j\wedge A_j+ \frac{\sqrt{-1}}{2} \cdot W^{pq} d\eta_p\wedge 
  d\bar{\eta}_q 
 \end{equation}
are compatible in the sense that $\Omega\wedge\bar{\Omega}=const\cdot\omega^{N}$.
\end{thm}

\begin{proof}
First we note that the local potential description of $V$ and $W$ by (\ref{eq:potential}) insures that

\begin{multline} \label{eq:complex_MA}
   dF_j= \frac{\sqrt{-1}}{2} \cdot d\left(\frac{\partial 
   W^{pq}}{\partial u_j} \right)\wedge d\eta_p\wedge d\bar{\eta}_q \\
   + \sqrt{-1} \cdot d\left( \frac{\partial V^{ij}}{\partial \eta_p} 
   \right)\wedge du_i\wedge d\eta_p   
   -\sqrt{-1} \cdot d\left( \frac{\partial V^{ij}}{\partial \bar{\eta}_q} 
   \right)\wedge du_i\wedge d\bar{\eta}_q \\
   = \frac{\sqrt{-1}}{2} \left( 
   \frac{\partial ^2W^{pq}}{\partial u_i\partial u_j}+
   4\frac{\partial ^2 V^{ij}}{\partial \eta_p \partial \bar{\eta}_q} \right) 
   du_i\wedge d\eta_p\wedge
   d\bar{\eta}_q =0.
\end{multline}
 Moreover, if $\theta_j$ denote the coordinates on the torus fiber such that $\partial/\partial \theta_j$ are the Hamiltonian vector fields, then the connection 1-forms  can be written up to exact forms on $B^\circ$ in terms of the local potential $\Phi$:
\begin{equation}
A_j=d\theta_j+\sqrt{-1}\left( \frac{\partial ^2 \Phi}{\partial u_j \partial \eta_p} d\eta_p -\frac{\partial ^2 \Phi}{\partial u_j \partial \bar{\eta}_q} d\bar{\eta}_q \right).
\end{equation}
And one can see explicitly that $F_j=dA_j$.

The integrability of the complex structure follows from the fact that the differential ideal generated by $(1,0)$-forms is closed:
\begin{multline}
  d(dz_j)=dV^{ij}\wedge du_i + \sqrt{-1}\cdot dA_j \\  
  =-\frac{\partial V^{ij}}{\partial u_k} du_i\wedge du_k
  -\frac{\partial V^{ij}}{\partial \eta_p} du_i\wedge d\eta_p
  -\frac{\partial V^{ij}}{\partial \bar{\eta}_q} du_i\wedge d\bar{\eta}_q \\
  - \left(\frac12 \frac{\partial W^{pq}}{\partial u_j}d\eta_p\wedge
  d\bar{\eta}_q +\frac{\partial V^{ij}}{\partial \eta_p}du_i\wedge d\eta_p 
  -\frac{\partial V^{ij}}{\partial \bar{\eta}_q}du_i\wedge d\bar{\eta}_q 
  \right)\\
  = \left( \frac12 \frac{\partial W^{pq}}{\partial u_j} d\bar{\eta}_q -2 \frac{\partial V^{ij}}{\partial \eta_p} du_i \right)\wedge d\eta_p,
\end{multline}
where we have only used 
 $\frac{\partial V^{ij}}{\partial u_k}=\frac{\partial V^{kj}}{\partial u_i}$.

It is equally easy to verify the K\"ahler condition: 
\begin{multline}
   d\omega=-du_j\wedge dA_j+ \frac{\sqrt{-1}}{2} \cdot dW^{pq}\wedge 
   d\eta_p\wedge  d\bar{\eta}_q \\
   = -\sqrt{-1}\left(\frac12 \frac{\partial W^{pq}}{\partial u_j}du_j\wedge 
   d\eta_p\wedge d\bar{\eta}_q \right) 
   +  \frac{\sqrt{-1}}{2}\frac{\partial W^{pq}}{\partial u_j}du_j\wedge 
   d\eta_p\wedge d\bar{\eta}_q = 0.
\end{multline}

Finally, the Ricci-flatness is manifest since $\det(h)=\det V^{-1} \cdot\det W =1$ in the complex coordinates $dz_j,d\eta_p$. 
\end{proof}

We are interested in applying the Gibbons-Hawking ansatz to description of the metrics on the toric Calabi-Yau hypersurfaces. Given a torus fibration of such hypersurface near the large complex structure point one may approximate the true Calabi-Yau metric by non-compact solutions of GH equations over different regions of the base. Away from the discriminant locus a semi-flat metric gives a good approximation. We try to argue that as the tori shrink to zero size the metric behavior near singular fibers can also be approximated by certain $y$-periodic solutions of GH ansatz. This local description is the main subject of the paper.

\subsection{Example: toric orbifold}
This is an important toy example which provides the local description of Gibbons-Hawking solutions for more interesting cases. Here for the standard toric orbifold metric one can actually write down an explicit solution to the Gibbons-Hawking equations. 

First we set up the notations. Let $N\cong\Z^{n+1}$ be an integral lattice in a real vector space $N_\R=N\otimes\R$. Denote by $N^*\subset N_\R^*$ the dual lattice in the dual space.

Let $\tau$ be an $n$-simplex with vertices $(w_0,w_1,\dots,w_n)$ in the lattice $N\cong\Z^{n+1}$, whose affine distance from the origin is 1. That is, there is a vector $\rho$ in the dual lattice $N^*$ such that $\<w_i,\rho\>=1$, all $i=0,\dots,n$.  Denote by $\mathcal T\subset N_\R$ the cone over $\tau$ and by $\mathcal T^\vee \subset N_\R^*$ the dual cone. 

Let $X_{\mathcal T}:=\operatorname{Spec}[z^m\suchthat m\in \mathcal T^\vee\cap N^*]$ be the associated affine toric variety (cf., e.g. \cite{Fu}). If $\Z\<w_0,\dots,w_n\>$ denotes the (finite index) sublattice in $N$ generated by $w_i$ and $G$ is the quotient group $N/\Z\<w_0,\dots,w_n\>$, then $X_{\mathcal T}$ is isomorphic to the orbifold $\C^{k+1}/G$. 

The real torus $N_\R/N$ acts on $X_{\mathcal T}$. But we will be interested rather in the action of its subtorus $T^n:=({N_\rho})_\R/{N_\rho}$, where $N_\rho:=\{v\in N\suchthat \<v,\rho\>=0\}$. The $n$-dimensional subspace $(N_\rho)_\R\subset N_\R$ can be naturally identified with the Lie algebra $\mathfrak t$ of $T^n$. The dual quotient space $N_\R^*/\rho$ is identified with ${\mathfrak t}^*$. 

Let $Q^\tau_i\subset N_\R^*/\rho $ denote the open normal cones to the vertices $w_i$ of $\tau$. Define a polyhedral complex $\Pi(\tau)$ in $N_\R^*/\rho$ to be the union of walls separating the $Q^\tau_i$'s:
$$\Pi(\tau):=\bigcup_{i\ne j}\mathrm{wall}_{ij},$$
with the orientation of each wall determined by the ordering of $\{i,j\}$.
Another way to look at $\Pi(\tau)$ is as being the image of the $(n-1)$-dimensional strata in $\mathcal T^\vee$ under the quotient map $N_\R^* \to N_\R^*/\rho$.

The vector $\rho$ lies in the interior of $\mathcal T^\vee$, hence $\eta=z^\rho:X_{\mathcal T}\to\C$ defines a regular function, which vanishes at the divisor in $X_{\mathcal T}$ corresponding to the boundary of ${\mathcal T}^\vee$. Together with the moment map $\mu:X_{\mathcal T}\to{\mathfrak t}^*$ this gives the torus fibration 
$$(\mu,\eta):X_{\mathcal T}\to N_\R^*/\rho \times \C,$$
whose restriction to $B^\circ := N_\R^*/\rho\times\C \setminus\nc(\tau)\times\{0\}$ is a principal $T^n$-bundle $\pi : M\to B^\circ$. 

To describe the topology of this bundle note that the homology group $H_2(B^\circ,\Z)\cong\Z^{n}$ is generated by the 2-spheres $\alpha_{ij}$ around every $\mathrm{wall}_{ij}$ oriented to match the orientation of $\mathrm{wall}_{ij}$. Looking at the relations among the generating cycles $\alpha_{ij}$ (abusing the notation) it is convenient to identify $H_2(B^\circ,\Z)$ with $\Lambda_\tau$, the (finite index) sublattice of $N_\rho$ generated by the elements $w_i-w_j$ for all pairs of $i,j$. On the other hand, the Chern classes of the restrictions of the $T^n$-bundle $M$ to the spheres $\alpha_{ij}$ are given by $w_i-w_j\in \mathfrak t$. Then the Chern class $[\nu]$  of the total bundle $M$:
$$[\nu]\in H^2 (B^\circ,\Z) \otimes N_\rho \cong \operatorname{Hom}(\Lambda_\tau,N_\rho)
$$ 
is the element given by the natural inclusion $\iota:\Lambda_\tau\hookrightarrow N_\rho$.

The final piece of notation before we describe the standard orbifold metric on $X_{\mathcal T}$ is the (finite index) sublattice $N'\subset N$ generated by $w_i$'s. Let $(N')^*\supset N^*$ be its dual lattice. Let $m_0,\dots,m_n$ be the minimal vectors in $(N')^*$ along the rays of ${\mathcal T}^\vee$. 

In polar coordinates the standard orbifold metric on the algebraic torus $(\C^*)^{n+1}\subset X_{\mathcal T}$ will be 
\begin{equation*}
h=\sum_{i=0}^n (d|z^{m_i}|+\sqrt{-1}|z^{m_i}| \<m_i,d\theta\>) \otimes (d|z^{m_i}|-\sqrt{-1}|z^{m_i}| \<m_i,d\theta\>).
\end{equation*}
The functions $z^{m_i}$ are defined only on the $|G|$-fold covering space of $X_{\mathcal T}$, but $|z|^{m_i}$ are well defined on $X_{\mathcal T}$ itself. So are the differential forms $\<m_i,d\theta\>$. 

To write this metric in the Gibbons-Hawking ansatz we choose a basis $\{e_i\}$ of $N_\rho= \mathfrak t_\Z$. Evaluating the moment map on the basis vectors  defines the coordinates $u_i=\mu(e_i)$ on $N_\R^*/\rho={\mathfrak t}^*$, thus giving an identification $N_\R^*/\rho$ with $\R^n$. The metric on each phase torus $T_{|z|}:=\{z\suchthat |z|^{m_i}=const\}\cong N_\R/N$ is constant, and, hence, it is given by a quadratic form $Q_{|z|}$ on $N_\R$. Let $(V^{-1})^{ij}$ be the matrix of restriction of $Q_{|z|}$ to $(N_\rho)_\R$ in the basis $\{e_i\}$. Then the functions $V^{ij}$ and $W:=\det V^{ij}$ give a solution to the GH equations.
 
Note that the top degree holomorphic form $\Omega_\tau$ coincide with the push forward  under the projection $\C^{n+1}\to\C^{n+1}/G$ of the standard volume form on $\C^{n+1}$. 

As an illustration free of orbifold complications let us write the ansatz for the standard Euclidean metric on $\C^{n+1}$ explicitly. In this case, $\tau$ is the standard $n$-simplex in $N:=\Z^{n+1}$, i.e. $w_i$ form a basis in $N$. We will fix the coordinates $z_i=|z_i|e^{i\theta_i}$ on $\C^{n+1}$.
The action of the torus  
$$T^n=\{(\theta_0,\dots,\theta_n) \suchthat \sum\theta_i=0 \}$$
on $\C^{n+1}$ gives rise to a principal $T^n$-bundle $M$ over $B^\circ=\R^n\times\C \setminus \nc(\tau)\times\{0\}$. Then, in the Gibbons-Hawking coordinates the metric on $M$ can be written as
$$
h=(V^{-1})^{ij}(V^{ki}du_k+\sqrt{-1}\cdot A_i)\otimes (V^{kj}du_k-\sqrt{-1}\cdot A_j) +W d\eta \otimes d\bar{\eta},
$$
where
\begin{gather*}
u_i=\frac12(|z_i|^2-|z_0|^2),\ i=1,\dots,n,\qquad \eta=z_0z_1\dots z_n,\\
W^{-1}=|z_0z_1\dots z_n|^2 \left(\frac{1}{|z_0|^2}+\frac{1}{|z_1|^2}+\dots+ \frac{1}{|z_n|^2}\right),\\
(V^{-1})^{ij}=|z_0|^2+\delta^{ij}|z_i|^2, \\
 A_j=d\theta_j-W|z_0z_1\dots \widehat{z_j}\dots z_n|^2 \cdot d(\theta_0+\theta_1+\dots+\theta_n).
\end{gather*}
The above expressions degenerate whenever two or more of the coordinates $z_i$ vanish. Thus, the discriminant locus $D\subset \R^n\times \C$ is given by $u\in\nc(\tau)$ and $\eta=0$. However, when written in the Euclidean coordinates the metric extends from $M$ to the standard flat metric $h=\sum_{i=0}^n dz_i \otimes d\bar z_i$ on $\C^{n+1}$. (To justify this it is enough to do the calculation for the associated K\"ahler form).

\subsection{Non-flat orbifold metrics}
We start with a $T^n$-torus bundle $\pi: M\to B^\circ$, where $B^\circ:=\R^n\times\C \setminus \nc(\tau)\times\{0\}$, which has the same topology as the orbifold bundle above. It is convenient to encode the topological information about the bundle by adding a distributional equation to the Gibbons-Hawking ansatz. Let $\gamma_\tau(u)$ be the $N_\rho$-valued 1-current in $\R^n$ supported on $\Pi(\tau)$ defined by
\begin{equation}
\gamma_\tau(\alpha)=\sum\limits_{s,s'} (w_s-w_{s'})\!\int\limits_{\mathrm{wall}_{ss'}}\!\alpha,
\end{equation}
for an $(n-1)$-form  $\alpha$. Once the basis of  $N_\rho$ is chosen we will write this current as $\gamma^j_\tau(u)$, where the superscript $j$ indicates its $N_\rho$-valuedness. Then replacing the equation (\ref{eq:complex_MA}) (which is automatic in $B^\circ$) for the closedness of the curvature forms by the distributional equation in $B=\R^n\times\C$:
\begin{equation} \label{eq:dist}
 \frac{\sqrt{-1}}{4\pi} \left( 
   \frac{\partial ^2 W}{\partial u_i\partial u_j}+
   4\frac{\partial ^2 V^{ij}}{\partial \eta \partial \bar{\eta}} \right)    du_i\wedge d\eta\wedge d\bar{\eta} =  \gamma^j_\tau(u)\wedge\delta(\eta)
\end{equation}
will guarantee that the fibration $\pi: M\to B^\circ$  has the right Chern class. This follows from integrating the curvature form $F$ over the generating 2-cycles $\alpha_{ij}\in H_2(B^\circ ,\Z)$:
$$ \int_{\alpha_{ij}} \frac1{2\pi} F= 
\int_{\beta_{ij}} \gamma_\tau(u)\wedge\delta(\eta) = w_i-w_j,
$$
where $\alpha_{ij}=\partial\beta_{ij}$ for a 3-chain $\beta_{ij}$ - a ball transversally intersecting $\mathrm{wall}_{ij}$ at a point. 

A remark on notation: here $\delta(\eta)$ stands for the two-current associated to the origin in $\C$ (the Dirac delta-function), but both $\gamma_\tau$ and $\delta$ in (\ref{eq:dist}) really mean the pullbacks of the corresponding currents to the product $\R^n\times \C$. We will continue to abuse this notation throughout the rest of the paper when there is no confusion possible.

\begin{lemma}\label{lemma:orbifold}
Suppose we have a Gibbons-Hawking solution on $\R^n\times\C \setminus \nc(\tau)\times\{0\}$, that is, a positive definite matrix function locally given by 
$$V^{ij}(u,\eta,\bar\eta)=\frac{\partial^2\Phi}{\partial u_j\partial u_j}
\quad \text{such that} \quad
W:=\det V^{ij} =-4\frac{\partial^2\Phi}{\partial \eta \partial \bar\eta},
$$ 
which satisfy the distributional equation (\ref{eq:dist}) in $\R^n\times\C$ and, in addition,
\begin{equation}\label{eq:complete}  
\int\limits_0^\infty \frac1{(V^{-1})^{ij} } du_i =\infty
\quad \text {for all $i$ and $j$}.
\end{equation}
Then the total space of the torus bundle $\pi: M\to \R^n\times\C \setminus \nc(\tau)\times\{0\}$ can be compactified to the fibration $\bar\pi: \bar M \to \R^n\times\C$ such that $\bar M$ is biholomorphic (in the orbifold sense) to $X_{\mathcal T}$ in a manner which respects the map
$$\eta:X_{\mathcal T}\to \C.$$ 
In particular, such solution defines a Ricci-flat K\"ahler metric on the orbifold $X_{\mathcal T}$ with the standard holomorphic volume form $\Omega_\tau$.
\end{lemma}

\begin{proof}
The discriminant locus $D=\nc(\tau)\times\{0\}$ is of codimension 3 in $\R^n\times\C$ and has a nice simplicial stratification. Using this stratification the topological compactification from $M$ to $\bar M=X_{\mathcal T}$ follows by extending the argument of \cite[Prop.~2.9]{Gr} to arbitrary dimensions and including the orbifold singularities. To prove matching of the complex structure we will follow closely \cite{LeB} where the argument is given for the $\C^2$ case. 

Let $\frac{\partial}{\partial \theta_j}$ be the Hamiltonian vector fields generating the $T^n$-action on $M$. Consider the commuting vector fields
\begin{equation}\label{eq:vf}
\xi_j=- \frac{\sqrt{-1}}{2}\left(\frac{\partial}{\partial \theta_j} -\sqrt{-1} J \left( \frac{\partial}{\partial \theta_j} \right) \right)= \frac12\left( (V^{-1})^{ij} \frac{\hat\partial}{\partial u_i} - \sqrt{-1}\frac{\partial}{\partial \theta_j}\right), 
\end{equation}
where 
$$ \frac{\hat\partial}{\partial u_i}= \frac{\partial}{\partial u_i} - A_k\left(\frac{\partial}{\partial u_i}\right)\frac{\partial}{\partial \theta_k},
$$
denotes the horizontal lift of $\frac{\partial}{\partial u_i}$, and $J: TM\to TM$ is the complex structure. Since the $\frac{\partial}{\partial \theta_j}$ preserve both the metric and the complex structure, it follows that the $\xi_j$ are holomorphic vector fields. Moreover, the flow of each $\xi_i$ is complete because of (\ref{eq:complete}). Hence, the  $\xi_j$ generate a holomorphic action of $(\C^*)^n$ on $\bar M$. 

The orbit structure of this action is easily seen to be identical with that of the toric $(\C^*)^n$-action on the orbifold $X_{\mathcal T}$. Namely, for each affine (real) codimension 2 plane
$$ L_a=\{(u,\eta)\suchthat \eta=a\} \subset \R^n\times\C,$$
the set $\pi^{-1}(L_a)$ is a union of orbits. For $a\ne 0$ the $\pi^{-1}(L_a)$ is a single orbit, and $\pi^{-1}(L_0)$ decomposes into $(n+1)$ orbits isomorphic to $(\C^*)^n$ plus a bunch of smaller dimensional orbits according to the polyhedral decomposition of $N^*_\R/\rho$ induced by $\Pi(\tau)$. 

Consider the following subsets of $N^*_\R/\rho\times\C$:
$$L_i:=\{(u,\eta)\suchthat \eta\ne 0\}\cup \{(u,\eta)\suchthat u\in Q^\tau_i\}.$$ 
We can cover the space $M$ by $(n+1)$ open sets $\mathcal U_i:= \pi^{-1}(L_i)$. Each map $\eta:\mathcal U_i\to \C$ defines a holomorphic principal $(\C^*)^n$-bundle over $\C$. Since $\C$ is Stein and contractible we conclude that each $\mathcal U_i$ is biholomorphic to the product $(\C^*)^n\times \C$.

Now notice that the bundle structures on $\mathcal U_i$ agree on their intersection $\pi: \bigcap_i \mathcal U_i\to\C^*$. Thus, $M$ is biholomorphic to the quotient space
$$\coprod_i \left( (\C^*)^n\times \C \right)/\sim,$$
where the equivalence relation on $(\C^*)^n\times \C^*$ is of the form
\begin{equation}
(z,\eta)_i \sim (f_{ij}(\eta)z,\eta)_j, \quad 0\le i,j \le n,
\end{equation}
for some cocycle $f_{ij}$ with values in the holomorphic maps $\C^*\to \operatorname{Aut} ((\C^*)^n)$. One can always choose the fiber coordinates $(z_k)_i$ in every $\mathcal U_i$ to be consistent with the vector fields $\{\xi_k\}$. That is, $\xi_k = (z_k)_i\frac{\partial}{\partial (z_k)_i}$ on the common part $\{\eta\ne 0\}\cong (\C^*)^{n+1}$ and $\eta$ is fixed by the $\xi_k$. In particular, this means that
$\int_{\gamma_k} d\log z_l = 2\pi\sqrt{-1}\, \delta_{kl}$,
where $\gamma_k$'s are the cycles generated by the $\xi_k$ flow.

Since the $(\C^*)^n$-action is the same on $(z)_i$ for different $i$ and  any automorphism of the principal homogeneous space $(\C^*)^n$ is given by an $n$-tuple of non-zero complex numbers we conclude that each $f_{ij}$ is an $n$-tuple of holomorphic functions $\C^*\to\C^*$. Another choice of trivializations $\mathcal U_i\cong (\C^*)^n\times \C$ amounts to modifying the cocycle $f_{ij}$ by a coboundary. That is, the biholomorphism type of $M$ is determined by the singularities of $f_{ij}(\eta)$ at $\eta=0$.

The $f_{ij}(\eta)$ cannot have essential singularities at 0, otherwise there would exist a sequence of points in $M$ converging to several distinct points, and, thus, $M$ would not be Hausdorff. On the other hand, the orders of vanishing of the $f_{ij}(\eta)$ can be read off from the residues of the 1-forms $d\log f_{ij}$, which is a topological information given by the Chern class of the bundle. More precisely,
$$\frac1{2\pi\sqrt{-1}}\oint d\log f_{ij} = [\nu](\alpha_{ij})=w_i-w_j.$$
Thus, the claimed biholomorphism is established on $M$, and it can be extended to $\bar M$ by an orbifold version of the Hartog's theorem.

In order to show that the holomorphic volume form is necessarily the standard one we will prove the following (stronger) statement. Given two sets of commuting holomorphic vector fields $\{\xi_i\}$ and $\{\zeta_i\}$ generating the $(\C^*)^n$-actions on $\eta:M\to\C$ which agree topologically, there is a biholomorphism $\phi: M\to M$ such that $\phi_* \xi_i=\zeta_i$. Then, taking $\zeta_i$ to be the standard torus action on $X_\mathcal T$ yields the claim about the volume form.

First, we will choose two sets of coordinates $\{\alpha_j,\eta\}$ and $\{\beta_j,\eta\}$ consistent with the topology of the action. Say, we take $(\alpha_j,\eta)$ to be $(z_j,\eta)_1$ with $\xi_j = \alpha_j\frac{\partial}{\partial \alpha_j}$ and $\eta$ fixed by the $\xi_j$, and similar for the $\{\beta_j,\eta\}$. Then the transition maps $f_{ij}(\eta)$ have the same zeros/poles for both sets of coordinates. Hence, the change-of-coordinates maps between $\{\alpha_j,\eta\}$ and $\{\beta_j,\eta\}$ extend from $(\C^*)^{n+1}$ to $M$ without zeros/poles. That is, if we write 
$$\alpha_j\frac{\partial}{\partial \alpha_j}=B_k^j \beta_k\frac{\partial}{\partial \beta_k},$$
then the matrix $B_k^j$ extends to an invertible matrix on $M$. 

We will be looking for a change of coordinates in the form $\beta_j=\alpha_je^{\psi_j}$ which would induce such a Jacobian matrix. Using the chain rule
\begin{equation*}
\alpha_j\frac{\partial}{\partial \alpha_j}=\left(\delta^j_k+\alpha_j\frac{\partial \psi_k}{\partial \alpha_j}\right) \beta_k\frac{\partial}{\partial \beta_k},
\end{equation*}
this amounts to solving the differential system
\begin{equation}\label{eq:diffsystem}
\alpha_j\frac{\partial \psi_k}{\partial \alpha_j}=B_k^j-\delta^j_k,
\end{equation}
which is, in general, overdetermined. 
However, in our case there are several restrictions on $B_k^j$. Namely, note that the forms 
$\frac{d\beta_j}{\beta_j}=B^k_j \frac{d\alpha_k}{\alpha_k}$ mod $d\eta$
are closed, that is 
\begin{equation}\label{eq:closed}
\alpha_i\frac{\partial B^k_j}{\partial \alpha_i}=\alpha_k\frac{\partial B^i_j}{\partial \alpha_k},
\end{equation}
 and have to satisfy the periodicity requirements of the action:
\begin{equation}\label{eq:periodic}
\frac1{2\pi\sqrt{-1}} \int_{\gamma_k} \frac{d\beta_j}{\beta_j}= \frac1{2\pi} \int_0^{2\pi} B^k_j (\alpha_1,\dots,\alpha_k e^{\sqrt{-1}\theta},\dots,\alpha_n,\eta)\, d\theta = \delta_j^k,
\end{equation}
where $\gamma_k$'s are the cycles generated by the $\alpha_k$ (or, equivalently, by the $\beta_k$) flow. Then, by writing $B^k_j$ in the power series form
\begin{equation}\label{eq:powerseries}
B^k_j=\sum\limits_{m_1 e_1+\dots+m_n e_n +r \rho \in \mathcal T^\vee} (b^k_j)_{m,r} \alpha_1^{m_1}\dots\alpha_n^{m_n}\eta^r, \quad m=\{m_1,\dots,m_n\} \in\Z^n,
\end{equation}
and using (\ref{eq:closed}) and (\ref{eq:periodic})
we conclude that $(b^k_j)_{m,r}$ have to satisfy the following restrictions:
\begin{itemize}
\item $m_i (b^k_j)_{m,r} = m_k (b^i_j)_{m,r}$;
\item $(b^k_j)_{m,r}=0$, if $m_k=0$ and $(m,r)\ne (0,0)$;
\item $(b^k_j)_{0,0}=\delta_j^k$.
\end{itemize} 
Thus, by setting $(a_j)_{0,r}:=0$ and $(a_j)_{m,r}:=\frac{(b^k_j)_{m,r}}{m_k}$ if there is an $m_k\ne 0$ (the ratio is independent of the choice of $k$) we can write down the power series solution to (\ref{eq:diffsystem}):
$$\psi_k=\sum\limits_{m_1 e_1+\dots+m_n e_n+r\rho\in \mathcal T^\vee} (a_k)_{m,r} \alpha_1^{m_1}\dots\alpha_n^{m_n}\eta^r,$$
that converges on the same polydisk  in $M$ as the power series (\ref{eq:powerseries}) for $B^k_j$ does. These functions $\psi_k$ provide the desired coordinate change.
\end{proof}

It is an interesting problem to exhibit the existence (and abundance) of the Gibbons-Hawking solutions. For instance, one can try to deform the standard (flat) orbifold metric. Applying continuity method techniques (cf. \cite[Ch.~17]{GTr}) this amounts to inverting a second order linear elliptic differential operator -- the linearization of the GH operator at the flat solution. The difficulty is that one of the eigenvalues of its symbol blows up at the discriminant.

In two dimensions (the original GH ansatz \cite{Ha},\cite{GH}) the equation becomes the usual Laplace equation. Since there are no positive harmonic functions on $\R^3$ except constants, any solution has the form 
$$V=\frac{\ell(\tau)}{2\sqrt{u^2+|\eta|^2}}+a,$$
where $\ell(\tau)$ is the length of $\tau$, and $a$ is a positive constant. This defines the famous Taub-NUT metric -- the first example of a non-trivial complete K\"ahler Ricci-flat metric on $\C^2$ and its quotients by cyclic groups.

\subsection{Periodic solutions}
The goal here is to set up the Gibbons-Hawking ansatz in such a way that the resulting complex manifold is identifiable with the local model for a Calabi-Yau toric hypersurface. There are no explicit solutions known in dimension higher than 2, unlike the orbifold case. But we will try to make use of the ansatz to get some information about the limiting behavior of solutions in certain degenerations.

We will adapt the notations from the orbifold example. Namely, $\tau=\{w_0,\dots,w_n\}$ is a simplex (but now of codimension $l+1$) in the lattice $N\cong\Z^{n+l+1}$, which has the distance 1 from the origin. Let $\sigma=\{v_0,\dots,v_l\}$ be a simplex in $N^*$ such that $\<\sigma,\tau\>=1$. In particular, it means that $\sigma$ also has distance 1 from the origin. Let $N_\sigma\subset N$ and $N^*_\tau\subset N^*$ be the sublattices orthogonal to $\sigma$ and $\tau$, respectively. And let $$N_\R^*/\sigma:=N_\R^*/\<v_0,\dots,v_l\>, \quad N_\R/\tau:=N_\R/\<w_0,\dots,w_n\>$$
be the corresponding (dual) quotients. 
The polyhedral complex $\Pi(\tau)$ provides a polyhedral decomposition of $N^*_\R/\sigma$ into cells $Q_i^\tau$ and $\gamma_\tau$ is the associated 1-current supported on $\Pi(\tau)$, as before. Also, we define the cone $\mathcal T:=\cone(\tau)$ in $N_\R$, its dual $\mathcal T^\vee\subset N^*_\R$ and let 
$$X_{\mathcal T}:=\operatorname{Spec}[z^m\suchthat m\in \mathcal T^\vee\cap N^*]$$
be the associated affine toric variety.

Every vertex $v_i\in\sigma$ lies in the interior of $\mathcal T^\vee$, and, hence the monomials $z^{v_i}$ belong to the coordinate ring of $X_{\mathcal T}$. Let $Z_{\sigma,\tau}$ denote the closure  of the affine hypersurface 
$$\{z\in (\C^*)^{n+l+1}\suchthat \sum_{i=0}^l z^{v_i} =1\}$$ 
in the toric variety $X_{\mathcal T}$.

The real torus $T^n:=(N_\sigma)_\R/N_\sigma$ acts on $X_{\mathcal T}$ and leaves the hypersurface $Z_{\sigma,\tau}$ invariant. We assume that this action is a holomorphic isometry and denote by $\mu:Z_{\sigma,\tau}\to N_\R^*/\sigma$ the corresponding moment map.

The natural inclusion $N^*_\tau\subset \mathcal T^\vee \cap N^*$ gives the projection $\kappa:X_{\mathcal T}\to T_\tau$ onto the algebraic torus $T_\tau:= \operatorname{Spec}[z^m \suchthat m\in N^*_\tau]$.
A choice of $\rho\in N^*$, such that $\<\rho,\tau\>=1$, defines the polynomial 
$$P_\sigma(z):=z^{-\rho}\sum_{i=0}^l z^{v_i},$$ 
which can be thought of as a polynomial in $T_\tau$. The zero divisor of $P_\sigma$ does not depend on the choice of $\rho$ and let $\Gamma_{\sigma}$ denote the 2-current in $T_\tau$ associated to $\{P_\sigma (z)=0\}$.

We will consider the map $(\mu,\kappa):Z_{\sigma,\tau}\to  N_\R^*/\sigma\times T_\tau$  as a torus fibration with the discriminant locus $D=\nc(\tau)\times \{P_\sigma (z)=0\}$.
When restricted to a domain $B \subset N_\R^*/\sigma\times T_\tau$ it defines a torus fibration $Z_{\sigma,\tau}(B) \to B$ which is a principal $T^n$-bundle
over $B^\circ:=B\setminus D$.  The Chern class is given, as before, by the inclusion $\iota:\Lambda_\tau\hookrightarrow N_\sigma$.

If the torus action is a holomorphic isometry a Ricci-flat metric on $Z_{\sigma,\tau}(B)$ can be written in the Gibbons-Hawking form. Our main goal of this section is to prove the converse. That is if we have a GH solution with the right Chern class then it defines a Ricci-flat metric on $Z_{\sigma,\tau}(B)$. 

To write everything in coordinates we choose a basis $\{e_i\}$ in $N_\sigma$ and a basis $\{m_p\}$ in $N^*_\tau$. This will defines the coordinates $u_i:=\mu(e_i)$ on $N_\R^*/\sigma\cong\R^n$ and $\eta_p:=\log(z^{m_p})=\<m_p,\log z\>$ on $T_\tau\cong (\C^*)^l$. 

\begin{defn}
Given a domain ${B}$ in $\R^n\times (\C^*)^l$ a {\it $(\sigma,\tau)$-type} solution to the Gibbons-Hawking ansatz in ${B}$ are two positive definite matrix functions -- a real $V^{ij}$ and a hermitian $W^{pq}$ -- on $B^\circ$ locally given by a potential:
$$V^{ij}=\frac{\partial^2\Phi}{\partial u_j\partial u_j}, \quad        
   W^{pq}=-4\frac{\partial^2\Phi}{\partial \eta_p \partial \bar{\eta}_q},
\qquad 1\le i,j \le n,\quad n+1\le p,q\le n+l, $$
such that $\det V^{ij}=\det W^{pq}$ and the distributional equation 
\begin{equation} \label{eq:distribution}
 \frac{\sqrt{-1}}{4\pi} \left( 
   \frac{\partial ^2W^{pq}}{\partial u_i\partial u_j}+
   4\frac{\partial ^2 V^{ij}}{\partial \eta_p \partial \bar{\eta}_q} \right)    du_i\wedge d\eta_p\wedge d\bar{\eta}_q = \gamma^j_\tau(u)\wedge\Gamma_{\sigma}(\eta)
\end{equation}
is satisfied in ${B}$. 
\end{defn}
The topological information about the bundle is again encoded in right hand side of the equation (\ref{eq:distribution}). 

To state the compatibility with the desired holomorphic volume form we recall  (cf., e.g., \cite{Bat}) that given an affine hypersurface $Z_f=\{f(z)=0\}\subset(\C^*)^{N+1}$ there is a distinguished nowhere vanishing top degree holomorphic form $\Omega_{CY}$ on $Z_f$, which is defined as a Poincar\'e residue of the meromorphic $(N+1)$-form $\frac{dz_0\wedge\dots\wedge dz_N}{f\cdot z_0 \dots z_N}$
on $(\C^*)^{N+1}$ with a single pole along $Z_f$. 

\begin{prop} Given a $(\sigma,\tau)$-type Gibbons-Hawking solution on a domain $B$, the total space of the torus bundle $M\to B^\circ$ can be compactified to the fibration $\bar M \to B$ such that $\bar M$ is biholomorphic (in the orbifold sense) to $Z_{\sigma,\tau}(B)$ in a manner which respects the fibration
$$\pi:Z_{\sigma,\tau}\to T_\tau.$$ 
In particular, such a solution defines a Ricci-flat K\"ahler (orbifold) metric on $Z_{\sigma,\tau}(B)$ with the holomorphic volume form $\Omega=\Omega_{CY}$.
\end{prop}
\begin{proof}
The topological compactification again can be drawn from \cite[Prop.~2.9]{Gr}. All we need to show is that a GH solution with asymptotics determined by (\ref{eq:distribution}) produces the right complex structure on $M$, which would then  uniquely extend to $\bar{M}$ by the orbifold version of Hartog's theorem. But this is a purely local question and it follows directly from Lemma~\ref{lemma:orbifold} dropping the completeness condition (\ref{eq:complete}) that becomes irrelevant.

To see matching of the volume forms let us choose local holomorphic coordinates $\{\tilde\eta_1,\dots,\tilde\eta_l\}$ on $(\C^*)^l$ such that the local equation for $\{P_\sigma=0\}$ is $\tilde\eta_1=0$. In these coordinates
the Gibbons-Hawking top degree holomorphic form on $M$ will be $\Omega=\Omega_\tau\wedge d\tilde\eta_2\wedge\dots \wedge\tilde\eta_l$, where $\Omega_\tau$ is the standard volume form on the $(n+1)$-dimensional toric orbifold $X_{\mathcal T}$ from Lemma \ref{lemma:orbifold}. Then, $\Omega$ is easily seen to coincide with the local expression for the distinguished form $\Omega_{CY}$ on $Z_{\sigma,\tau}$.
\end{proof}

\section{Limiting behavior of solutions}
\subsection{Exponential decay lemma}
We would like to analyze the behavior of Gibbons-Hawking solutions when the tori (both in the fibers and in the base) are shrinking. First, let us introduce a non-linear differential equation of the Monge-Amp\`ere type.

\begin{defn}
We refer to a pair of real positive definite matrix functions $V^{ij}, W^{pq}$ as a solution to the {\it split} Monge-Amp\`ere equation in an open subset $R\subset \R^n\times\R^{l}$ if $\det V^{ij}=\det W^{pq}$ and they are locally given by a smooth potential function $K$:
\begin{equation}\label{eq:splitMA}
V^{ij}=\frac{\partial^2 K}{\partial s_i\partial s_j},\quad W^{pq}=-\frac{\partial^2 K}{\partial t_p \partial t_q},
\quad 1\leq i,j \leq n,\quad n+1\leq p,q \leq n+l.
\end{equation}
\end{defn}

To describe the asymptotics at the discriminant we would like to treat the simplex $\sigma\subset N^*$ on the same footing as $\tau$. Namely, we let $\Sigma\subset N^*_\R$ be the cone over $\sigma$, and let $\Sigma^\vee$ be its dual cone in $N_\R$. The polyhedral complex $\Pi(\sigma)$ provides a polyhedral decomposition of $N_\R/\tau$ into cells $Q_i^\sigma$. Denote by  $\gamma_\sigma$ the $N^*_\tau$-valued 1-current defined in the same way as $\gamma_\tau$.

\begin{defn} Given a domain $R$ in $N^*_\R/\sigma\times N_\R/\tau\cong\R^n\times\R^l$ a $(\sigma,\tau)$-type {\em singular} solution to the split Monge-Amp\`ere equation in $R$ is a pair of matrix functions $V^{ij}, W^{pq}$ which are local Monge-Amp\`ere solutions in $R\setminus (\nc(\tau)\times\nc(\sigma))$ with asymptotics at the discriminant locus governed by the distributional equation
\begin{equation} \label{eq:splitMAdistribution}
 \left( 
   \frac{\partial ^2W^{pq}}{\partial s_i\partial s_j}+
   \frac{\partial ^2 V^{ij}}{\partial t_p \partial t_q} \right) 
   ds_i\wedge dt_p= \gamma^j_\tau(s) \wedge \gamma^q_\sigma(t).
\end{equation}
\end{defn}

\begin{conjecture}[Exponential decay lemma]\label{conj:decay}
Given a convex domain  $R$ in $\R^k\times\R^l$ and a $(\sigma,\tau)$-type solution $V, W$ of the split Monge-Amp\`ere equation in $R$ there is a real one-parameter family of $(\sigma,\tau)$-solutions $V_\lambda, W_\lambda$ to the Gibbons-Hawking ansatz in $\lambda R\times(S^1)^l$ such that
\begin{itemize}
\item The diameter of the circles both in the fiber $T^n$ and in the torus part $(S^1)^l$ of the base away from the discriminant is roughly given by $\lambda^{-1}$.
\item The zero Fourier modes of the GH solutions $V_{\lambda}^0 (u,x), W_{\lambda}^0 (u,x)$ as functions of the rescaled variables $s,t$, where $u=\lambda s, x=\lambda t$, will converge (in some properly weighted norm on the function space) to $V(s,t), W(s,t)$ as $\lambda\to\infty$.
\item The higher Fourier modes decay exponentially away from the discriminant $\nc(\tau)\times\nc(\sigma)$ in $\lambda R$, uniformly in $\lambda$. That is, if $\beta(u,x)$ denotes the Euclidean distance from the point $(u,x)$ to the discriminant, then
$$|V_{\lambda}^m(u,x)| \le C_1e^{-\beta(u,x)|m|}, \quad |W_{\lambda}^m(u,x)| \le C_2e^{-\beta(u,x) |m|}, $$
for some constants $C_1,C_2$, and large enough $\lambda$ and $\beta$.
\end{itemize}
\end{conjecture}

We would like to give some simple examples and a rough argument based on those why we believe this conjecture is true. Note, however, that once justified, it will have an important consequence for the metric collapse program for the toric hypersurfaces and complete intersections: 
\begin{cor}
The metric space $(Z_{\sigma,\tau}(\lambda^{-1}R),\lambda^{-2} g_\lambda)$, where $g_\lambda$ is the Riemannian (orbifold) metric from the Gibbons-Hawking ansatz, converges in the Gromov-Hausdorff sense to $(R,g^{ij}_\infty)$, with the limiting metric $g^{ij}_\infty= V^{ij}ds_ids_j+W^{pq}dt_pdt_q$.
\end{cor}

\subsection{The semi-flat case}
We consider the case when either $l=0$, or $n=0$. In both situations the discriminant locus is empty and the total space $M$ is just the product of the domain $R$ and the torus $T^n$. We can use any solution of the classical real Monge-Amp\`ere equation  in $R$ and extend it to a Gibbons-Hawking solution on $M$ by setting higher Fourier modes to zero. 
In the obvious complex structure  this will give a Ricci-flat metric on $M=\bar{M}$ (cf. \cite{Hi}, \cite{Leung}, \cite {LYZ}).

\subsection{Two dimensional example: local K3 (after \cite{OV} and \cite{GW})}
This is the periodic version of the original Gibbons-Hawking ansatz \cite{GH},\cite{Ha}. We consider the case when $n=l=1$ and both simplices $\tau$ and $\sigma$ are of length 1, although the construction works for a non-unimodular case as well. 

The Gibbons-Hawking equation in this case is equivalent to the Laplace equation for $V(u,x,y)$ ($=W(u,x,y)$) on a domain in the cylinder $\R\times\R\times S^1$ with the Dirac $\delta$-function on the right hand side. We can write both the solution $V(u,x,y)$ and the $\delta$-function in the Fourier expansion:
$$ V(u,x,y)=\sum_{m\in\Z} V_m e^{i my}, \quad \gamma^j_\tau(u)\delta_{P_\sigma}(\eta)=-\delta(u,x,y)=- \frac1{2\pi} \sum_{m\in\Z} \delta(u,x) e^{i my}.
$$
Here the minus sign takes into account the orientation of the circle action when passing from currents to generalized functions. 

The Gibbons-Hawking equation
$$ \frac{\partial^2 V_\lambda}{\partial u^2}+ \frac{\partial^2 V_\lambda}{\partial x^2}+\frac{\partial^2 V_\lambda}{\partial y^2}=-2\pi \cdot \delta(u,x,y)$$
decomposes into the Helmholtz equations according to the Fourier modes:
$$ \frac{\partial^2 V_{\lambda}^m}{\partial u^2}+ \frac{\partial^2 V_{\lambda}^m}{\partial x^2}- m^2 V_{\lambda}^m= - \delta(u,x), \quad m\in\Z.$$

On the other hand, the $(\sigma,\tau)$-type split Monge-Amp\`ere equation in the rescaled coordinates $s=\lambda^{-1}u, t=\lambda^{-1} u$ is the two-dimensional Laplace equation:
$$ \frac{\partial^2 V}{\partial s^2 }+ \frac{\partial ^2 V}{\partial t^2}=- \delta(s,t),$$
whose fundamental solutions are in the form $V(s,t)=-\frac1{4\pi} \log|s^2+t^2|+h(s,t)$, for harmonic $h$. Thus, one can take the zero mode of the corresponding Gibbons-Hawking solution to be $V_{\lambda}^0 (u,x)=V(\lambda^{-1} u,\lambda^{-1} x)$, as long as $V(s,t)$ stays positive in $R$.  
As for the higher modes, it is known that a fundamental solution to the Helmholtz equation with $m\ne 0$ may be given by the Bessel function
$$V_\lambda^m= K_0(|m| r) \sim {\sqrt{\frac\pi{2|m|r}}} \,e^{-|m|r}\left(1+O(r^{-1})\right), \text{ where } r^2=u^2+x^2,
$$
which decays exponentially as required.

\subsection{Higher dimensional case}
The full proof of the conjecture in this general case will probably require some very non-trivial application of the continuity method to deform the given split solution, then to introduce exponentially small higher modes and do some clever estimates afterwords. Meanwhile, we want to indicate a rough argument why some of the ideas from the K3 example above may still work in general.

To have the Fourier modes of the solutions defined on the same domain, independent of $\lambda$, we can scale the variables by $\lambda$:
$$s_i=\frac{u_i}{\lambda}, \quad t_p=\frac{x_p}{\lambda}, \quad \tilde{y_p}=\frac{y_p}{\lambda}.$$
Then the GH solutions are on $R\times (S^1/\lambda)^l$ for all $\lambda$, and their Fourier modes are functions on $R$. To keep up with the complex structure the logarithmic map $(\C^*)^l \to \R^l$ has to scale by $\lambda$ as well: 
$$\log_{e^\lambda} (z_1,\dots, z_n) := \frac1\lambda (\log|z_1|,\dots,\log|z_n|).$$

We would like to recall a few basic facts from ``tropical'' geometry (cf., e.g., \cite{Mi}). Given a polynomial $P_\sigma(z)$ in $(\C^*)^l$ the amoeba $\mathcal A^\lambda_\sigma$ is defined to be the image of the rescaled log map:
$$\mathcal A^\lambda_\sigma:= {\log_{e^\lambda} (\{P_\sigma=0)\}}\subset \R^l.$$
As $\lambda\to\infty$ the amoeba $\mathcal A^\lambda_\sigma$ approaches its spine $\mathcal A^\infty_\sigma=\Pi(\sigma)$.
The Ronkin function 
$$N_{\sigma}(x):=\frac1{(2\pi\sqrt{-1})^l}\int\limits_{\log|z|=x}\! \log|P_\sigma(z)|^2 \ \frac{dz_1}{z_1}\wedge \dots \wedge\frac{dz_l}{z_l}$$
is defined up to a linear function, which depends on a particular choice of $\rho\in N^*$ used in the definition of the polynomial $P_\sigma$. Denote by 
$$N_\sigma^\lambda(t):=\frac1\lambda N_\sigma(\lambda t)$$
the rescaled Ronkin function. The point is that $N^\lambda_{\sigma}(t)$ is a continuous function, linear on each connected component of $\R^l\setminus A^\lambda_\sigma$, with slopes given by the $v_i$. As $\lambda\to \infty$, it converges to the piece-wise linear function $N_\sigma^\infty(t)$ whose corner locus is $\Pi(\sigma)$ with the slopes $v_i$ over the $Q^\sigma_i$. In particular, the Hessian of $N_\sigma^\infty(t)$ can be interpreted as the generalized function associated with the current $\gamma_\sigma$.

We would like to analyze the right hand side of the equation (\ref{eq:distribution}) written in the Fourier expansion. The factor $\gamma_\tau$ carries over to every mode, while 
$$\Gamma_\sigma=\frac{\sqrt{-1}}{2\pi}\partial\bar\partial \log|P_\sigma|^2$$
decomposes into currents supported on $\mathcal A^\lambda_\sigma$. In particular, since the exterior differentiation commutes with averaging, we conclude that the zero mode of $\Gamma_\sigma$ is given by the Hessian of the Ronkin function:
\begin{multline*}
\Gamma_\sigma^0=\frac1{(2\pi\sqrt{-1})^l}\!\int\limits_{\log|z|=x}\! \Gamma_\sigma \ \frac{dz_1}{z_1} \wedge \dots \wedge \frac{dz_l}{z_l}=\\ 
\frac{\sqrt{-1}}{2\pi} \cdot \frac1{(2\pi\sqrt{-1})^l}\  \partial\bar\partial  \! \int\limits_{\log|z|=x}\! \log|P_\sigma|^2 \ \frac{dz_1}{z_1}\dots  \frac{dz_l}{z_l}= \frac1{4\pi}
\frac{\partial^2 N_{\sigma}(x)}{\partial x_p \partial x_q} dx_p\wedge dy_q,
\end{multline*}
where $\frac{dz_q}{z_q}=dx_q+\sqrt{-1} dy_q$. Then substituting $x=\lambda t$ yields
$$\frac{\partial^2 N_{\sigma}(x)}{\partial x_p \partial x_q} dx_p=  \frac{\partial^2 N^\lambda_{\sigma}(t)}{\partial t_p \partial t_q} dt_p.
$$
But, as $\lambda\to \infty$, the Ronkin function $N^\lambda_{\sigma}(t)$ tends to its tropical limit $N_\sigma^\infty(t)$, and the same is true for the respective Hessians. With the distributional interpretation of $\Hess N_\sigma^\infty(t)$ this means that the current $\Gamma_\sigma^0$ converges to $\gamma_\sigma dy_q$.

As for the higher modes, we note that since $\tilde{y}$ is now $(2\pi\lambda^{-1})$-periodic, there is a factor of $\lambda^2$ in the zero order term of the Helmholtz-type equation for $m\ne 0$. By analogy with the Bessel functions we hope that the spectral theory will force the higher modes decay exponentially away from the locus $\mathcal A^\lambda_\sigma \times\Pi(\tau)$ with the exponent roughly proportional to $\lambda$.

One can go about proving the exponential decay conjecture by starting with the given split Monge-Amp\`ere solution and constructing a family of solutions but with the factor $\gamma_\sigma dy_q$ in the right hand side being replaced by a more regular $\Gamma_\sigma^0$. This will give a family of semi-flat Gibbons-Hawking solutions. Then one can argue that since the higher modes can be taken exponentially small, they may be considered, in some sense, as perturbation of the semi-flat solution.

\section{Local mirror symmetry and Legendre transform}
\subsection{Linear algebra of Legendre transform and Monge-Amp\'ere equations}
The classical fact, implicitly used in \cite{GW}, is that solving the two-dimensional Laplace equation is equivalent  to solving the real two-dimensional Monge-Amp\`ere equation. This can be easily generalized to higher dimensions. Namely, the chain rule and elementary linear algebra for a particular coordinate transformation yields the following. 

\begin{lemma} \label{lemma:Legendre}
If $K(s,t)$ is a local solution
to the split Monge-Amp\`ere equation (\ref{eq:splitMA}), then $\Psi(y)$ is a (local) solution of the classical (real) Monge-Amp\`ere
equation  
\begin{equation} 
 \det \frac{\partial^2\Psi}{\partial y_i \partial y_p}=1,
 \quad 1\leq i,j \leq n+l, 
\end{equation}  
where 
\begin{equation}\label{eq:coord}
y_i=\frac{\partial K}{\partial
s_i},\ 1\leq i\leq n,\qquad y_p=t_p,\ n+1\leq p\leq n+l,
\end{equation}
and $\Psi(y)$ is the partial
Legendre transform of $K(s,t)$, defined by 
\begin{equation}
\frac{\partial \Psi}{\partial
y_i}=s_i,\ 1\leq i\leq n, \qquad \frac{\partial \Psi}{\partial y_p}=-\frac{\partial K}{\partial t_p}, \ n+1\leq p\leq n+l.
\end{equation}
\end{lemma}

\begin{proof}
First, we check the very existence of the partial Legendre transform. Consider the following Jacobian and Hessian matrices:
\begin{equation}
\frac{\partial (y_i,y_p)}{\partial (s,t)}
=\left(
\begin{array}{cc}
\frac{\partial^2 K}{\partial s_j \partial s_i} & \frac{\partial^2 K}{\partial t_q \partial s_i}\\
0 & \mathbbm{1}
\end{array}\right)
=\left( \begin{array}{cc} V & B\\ 0 & \mathbbm{1} \end{array}\right), \quad 
\Hess K (s,t) =\left( \begin{array}{cc} V & B\\ ^tB & -W \end{array}\right).
\end{equation}
Then
\begin{equation} 
\frac{\partial (s,t)}{\partial (y_i,y_p)}
=\left( \begin{array}{cc} V^{-1} & -V^{-1}B\\ 0 & \mathbbm{1} \end{array}\right),  \quad 
\Hess\Psi(y) =\left( \begin{array}{cc} 
V^{-1} & -V^{-1}B \\ ^t(-V^{-1}B) & W + {^tB} V^{-1}B \end{array}\right). 
\end{equation}
This shows that if both $V$ and $W$ are symmetric and positive definite, then $\Hess\Psi$ is also a symmetric positive definite matrix, so there exists (locally) a convex function $\Psi$ -- the partial Legendre transform.

On the other hand the inverse of a non-degenerate $2\times 2$ block matrix 
 with $\det A \ne 0$ and $\det D \ne 0$ is:
\begin{equation}
\left( \begin{array}{cc} A & B\\ C & D \end{array}\right)^{-1} 
=\left( \begin{array}{cc} (A-BD^{-1}C)^{-1} & (-A+BD^{-1}C)^{-1} BD^{-1}\\ 
(-D+CA^{-1}B)^{-1} CA^{-1} & (D-CA^{-1}B)^{-1} \end{array}\right),
\end{equation}
which can be checked directly  by multiplying the right-hand side by the matrix
$\left( \begin{smallmatrix} A & B\\ C & D  \end{smallmatrix} \right) $. 
Another useful linear algebra fact is that for any invertible matrix $M$ any minor of $M^{-1}$ is equal (with the usual parity sign) to the complementary minor of $M$ divided by $\det M$.
This implies that 
\begin{equation}
\det \left( \begin{array}{cc} A & B\\ C & D \end{array}\right) =1 \quad \Longleftrightarrow \quad \det A^{-1}=\det (D-CA^{-1}B).
\end{equation} 
Applied to $\Hess \Psi$ the last observation shows that $\Psi$ is a local Monge-Amp\`ere solution if and only if $\det V= \det W$.
\end{proof}

\subsection{Monge-Amp\`ere manifolds}
\begin{defn}
A Riemannian manifold $(Y,g)$ is called {\it Monge-Amp\`ere} if it possesses an integral affine structure and the metric is, locally in affine coordinates, given by a potential $\Psi$, that is, $g_{ij}= \frac{\partial^2\Psi}{\partial y_i \partial y_j}$, that satisfies the real Monge-Amp\`ere equation  $\det g_{ij}=1$. 
\end{defn}

Cheng and Yau \cite{ChY} proved that every {\em compact} Monge-Amp\`ere manifold is diffeomorphic to a torus and the Monge-Amp\`ere structure is a deformation of the standard flat structure on $\R^n/\Z^n$. To enrich this fairly boring class of manifolds we will allow certain singularities along a codimension 2 discriminant locus. 

The Monge-Amp\`ere condition in the affine geometry is an analog of the Ricci-flatness condition on a K\"ahler manifold with similar global rigidity properties. Roughly speaking, one should expect some sort of the Calabi conjecture saying that each topological class of the metric has a unique Monge-Amp\`ere representative (cf. \cite{HZh2} and \cite{KTch}).

Here we are concerned with the local picture in which many Monge-Amp\`ere structures may exist. However, for applications to the metric collapse of toric Calabi-Yau hypersurfaces we will be interested only in a certain special class of the Monge-Amp\`ere structures, called bi-PIKAS in \cite{HZh2}.  We will show that the latter always arise from singular split Monge-Amp\`ere solutions of the corresponding type.

The base $(\R^n\times\R^l)\setminus (\nc(\tau)\times\nc(\sigma))$ has an open covering $\{U_{v_i}, U_{w_j}\}$ where $U_{v_i}:= \R^n\times Q^\sigma_i$ and $U_{w_j}:= Q^\tau_j\times \R^l$. The nerve of this covering is the complete bipartite graph on the vertices $v_i$ of $\sigma$ and $w_j$ of $\tau$. We assume that the affine structure on $(\R^n\times\R^l)\setminus (\nc(\tau)\times\nc(\sigma))$ is polyhedral of type $(\{U_v\},\partial\Sigma^\vee)$. That is,  there is a homeomorphism $\R^n\times\R^l\to\partial\Sigma^\vee$ which provides affine coordinates on every chart $U_{v_i}$ by identifying it with the maximal dimensional face of $\partial\Sigma^\vee$ orthogonal to $v_i$. Also we assume that the {\em dual} affine structure is polyhedral of type $(\{U_w\},\partial\mathcal T^\vee)$. Moreover, the transformation maps between the affine coordinates on $U_{v_i}$ and $U_{w_j}$ are assumed to be the natural projections from the subspaces $v_{i}^\perp := \{y \in N_\R \suchthat \<v_{i},y\> = 0 \}$ to the quotients $N_\R/w_j$. This data constitutes a bi-polyhedral integral K\"ahler affine structure (bi-PIKAS for short) on $(\R^n\times\R^l)\setminus (\nc(\tau)\times\nc(\sigma))$ of type $(\partial\Sigma^\vee, \partial\mathcal T^\vee)$ (cf. \cite{HZh2}). We will abbreviate the type to be $(\sigma, \tau)$.

It is straight forward to see that the monodromy of the (linear part of the) affine structure along a loop $(v_{i_1}w_{j_1} v_{i_2}w_{j_2})$ is given by $y \mapsto y + \<v_{i_2}-v_{i_1},y\>(w_{j_2}-w_{j_1})$ (cf. \cite{HZh} for details). Since the radiance obstruction class vanishes locally (cf. \cite{GS}) we can always choose the affine coordinates such that there is no translational part of the monodromy.

Given a convex domain $R$ in $\R^n\times\R^l$ which contains the origin, a Monge-Amp\`ere bi-PIKAS on $R$ of type $(\sigma, \tau)$ will be a restriction of a Monge-Amp\`ere bi-PIKAS of the same type on $\R^n\times\R^l$.

\begin{thm}
There is a bijection between the sets of $(\sigma,\tau)$-type split Monge-Amp\`ere solutions on a domain $R$ and $(\sigma,\tau)$-type Monge-Amp\`ere bi-PIKAS on $R$. The Riemannian metric on $R^\circ:=R\setminus \Pi(\sigma)\times \Pi(\tau)$ is given by $V^{ij} ds_i ds_j + W^{pq} dt_p dt_q$, where $(V^{ij},W^{pq})$ is the corresponding split MA solution.
\end{thm}
\begin{proof}
Let $K_\alpha$ be local potentials for a given split Monge-Amp\`ere solution $(V,W)$. We will use the partial Legendre transform to define new coordinates $y_i,y_p$ as in Lemma~\ref{lemma:Legendre}. We claim that these are {\em affine} coordinates, that is the transition maps are affine linear. 

Indeed, by comparing the two local potentials in the overlap $U_\alpha\cap U_\beta$ we see that $K_\alpha-K_\beta$ have to be affine linear in both $s$ and $t$ variables. Thus $y_\alpha-y_\beta$ are affine functions of $t_p$, hence affine functions of $y_p$.  

Lemma~\ref{lemma:Legendre} also guarantees that $\Psi_\alpha$ are local Monge-Amp\`ere potentials in the affine coordinates $y$. The polyhedral property follows from noticing that the charts $U_v$ are bounded by linear inequalities in $t_p$, and hence by (the same) linear inequalities in the new (affine) coordinates $y_p$. 

Applying the partial Legendre transform to the other half of the variables $(s,t)$ gives the {\em dual} MA structure on $R^\circ$. Same considerations as above show that the dual affine structure is polyhedral of type $(\{U_w\},\partial\mathcal T^\vee )$.

The final step is to show that the resulting Monge-Amp\`ere bi-PIKAS has the right type, that is, to compute its monodromy around a loop $(v_{i_1}w_{j_1} v_{i_2}w_{j_2})$. For this we consider the closed $(N^*_\tau)_\R\otimes (N_\sigma)_\R$-valued 1-form on $R^\circ$:
$$\beta ^{iq}= \frac{\partial W^{pq}}{\partial s_i} dt_p-
   \frac{\partial V^{ij}}{\partial t_q} ds_j.$$
Making use of the distributional equation (\ref{eq:splitMAdistribution}) and applying Stokes' theorem we can compute the 
integral of $\beta$ along the loop $(v_{i_1}w_{j_1} v_{i_2}w_{j_2})=:\partial S$
\begin{equation}\label{eq:stokes}
\oint_{\partial S} \beta= \iint_S d\beta = \iint_S \gamma_\sigma \gamma_\tau 
=-(v_{i_2}-v_{i_1})\otimes(w_{j_2}-w_{j_1}),
\end{equation} 
where the minus sign is to take care of orientation of the loop.
Thus, $\beta$ has holonomy which depends only on the homotopy class of the path. That is, $\beta^{iq}=d f^{iq}$ for a multi-valued function $f^{iq}=-\frac{\partial^2 K}{\partial s_i \partial t_q}$. 

The affine coordinates $y_p,\ n+1\le p\le n+l,$ are single valued. Hence, there is no monodromy in $dy_p$. On the other hand, the ambiguity in the remaining differentials
$$dy_i=\frac{\partial^2 K}{\partial s_i \partial t_q} dt_q + \frac{\partial^2 K}{\partial s_i \partial s_j} ds_j=V^{ij} ds_j - f^{iq} dt_q$$
is coming exactly from the multi-valuedness of $f^{iq}$. Thus, the monodromy around the loop $(v_{i_1}w_{j_1} v_{i_2}w_{j_2})$ is given by $\mathbbm{1} + (v_{i_2}-v_{i_1})\otimes(w_{j_2}-w_{j_1})$.

Tracing backwards through the above argument shows that the converse is true.  Namely, given a Monge-Amp\`ere bi-PIKAS of type $(\sigma,\tau)$ we can use it together with its dual to define single-valued coordinates $(s,t)$ on $R^\circ$ which will obviously extend to $R$. Moreover, due to the polyhedral properties of these bi-PIKAS, the discriminant locus (and, thus, the support of the current in the right-hand side of (\ref{eq:splitMAdistribution})) in $R$ will be exactly given by $\dU_v\cap\dU_w=\nc(\sigma)\times\nc(\tau)$. Applying the Stokes' formula (\ref{eq:stokes}) with the prescribed monodromy of the affine structure to every 2-cycle in $R\setminus \nc(\sigma)\times\nc(\tau)$ identifies the current as $\gamma_\sigma \gamma_\tau $ and, hence, guarantees the $(\sigma,\tau)$-type asymptotics of the Monge-Amp\`ere potential $K(s,t)$. 
\end{proof}

\begin{remark}
The discriminant locus $D=\dU_v\cap\dU_w$ is not expected to be affine linear. Even though $D$ lies in the polyhedral boundary $\dU_v$, the boundary $\dU_w$ will be wiggled in the affine coordinates unless the partial Legendre transform is linear.
\end{remark}

The proof of the above theorem shows, in fact, that the distributional equation (\ref{eq:splitMAdistribution}) guarantees the potential form of the metric. One just needs $V$ and $W$ to be symmetric. Indeed, locally we can define closed 1-forms $dx_q:=f^{iq}ds_i+W^{pq} dt_p$ (the $(s,x)$ will form the dual affine coordinates). Then the form $dK:= y_i ds_i + x_q dt_q$ is also closed, and $K$ can be used as local potential: $V^{ij}=\frac{\partial^2 K}{\partial s_j\partial s_j}$ and $W^{pq}=-\frac{\partial^2 K}{\partial t_p \partial t_q}$.

We would like to finish by mentioning an obvious application of the construction in this section to mirror symmetry. As was noted in \cite{Hi} the mirror duality in the semi-flat case is provided by the Legendre transform. This statement continues to hold in the neighborhood of the discriminant locus as well. Namely, in the single-valued coordinates (which are not affine) the full Legendre transform takes a singular split Monge-Amp\`ere solution of type $(\sigma,\tau)$ into that of type $(\tau,\sigma)$.    

\bibliographystyle{alpha}
\bibliography{cy}

\end{document}